\documentclass[11pt,a4paper]{article}
\usepackage[english]{babel}
\usepackage{fullpage}
\usepackage{indentfirst}
\usepackage{amsmath,amsfonts,amsthm}
\usepackage[colorlinks,linkcolor=blue]{hyperref}
\usepackage[abbrev,nobysame]{amsrefs}
\usepackage{tikz}
\usetikzlibrary{decorations.markings,intersections, calc, backgrounds}
 \usepackage{xfrac}
\usepackage{mathtext}
\usepackage{moresize}
\usepackage{pdfpages}
\usepackage{bm}
\usepackage{multicol}
\usepackage{wrapfig}
\usepackage{graphicx}
\usepackage{subcaption}
\usepackage{addfont}
\usepackage{mathtools}
\usepackage{amssymb}

\newtheorem{theorem}{Theorem}[section]
\newtheorem{lemma}{Lemma}[section]
\newtheorem{cor}[lemma]{Corollary}

\newtheorem{example}{Example}[section]
\theoremstyle{definition}
\newtheorem{definition}[lemma]{Definition}
\theoremstyle{remark}
\newtheorem{remark}[lemma]{Remark}

\newcommand{\R}{\mathbb R}

\newcommand{\ga}{\gamma}
\newcommand{\pd}{\partial}
\newcommand{\co}{\colon}
\newcommand{\LL}{\mathcal{L}}
\newcommand{\ov}{\overline}
\newcommand{\Norm}{{\|{\cdot}\|}}
\newcommand{\8}{^*}
\newcommand\be{\begin{equation}}
\newcommand\ee{\end{equation}}
\numberwithin{equation}{section}

\title{On the Finiteness of Reflections in Minkowski Billiards in a Polyhedral Angle}
\date{}
\author{R.~V.~Barinov$^1$}
\begin{document}
\maketitle
\footnotetext[1]{Saint Petersburg University, 7/9 Universitetskaya nab., St. Petersburg, 199034 Russia}
\renewcommand{\thefootnote}{}
\footnotetext{This work was performed at the Saint Petersburg Leonhard Euler International Mathematical Institute and supported by the Ministry of Science and Higher Education of the Russian Federation (agreement no. 075–15–2025–343).}
\begin{abstract}
    We prove that in a Minkowski billiard in a polyhedral cone, that is a billiard in a polyhedral cone in a normed space with a smooth strictly convex norm, any billiard trajectory has only finitely many reflections off the walls. We also give an example of a billiard that admits a trajectory with arbitrarily many reflections.
\end{abstract}
\section{Introduction}
\qquad In the 1970s, Ya.~G.~Sinai and G.~A.~Galperin formulated the following problem, motivated by the study of ideal gas models: is the number of collisions bounded in a system of $n$ hard balls in Euclidean space $\R^d$ where the balls move uniformly and collide elastically? The set of admissible positions of the balls is naturally identified with the space $\R^{dn}$ with $\frac{n(n-1)}2$ cylindrical sets removed, whose points correspond to configurations where some pairs of balls overlap. Upon an appropriate linear transformation, the evolution of the system becomes a billiard trajectory in this space, with the standard law ``the angle of incidence equals the angle of reflection''; see, e.g., \cite{Gal81}.

For a one-dimensional ideal gas, that is a system of $n$ point masses moving uniformly along a line and colliding elastically, the problem reduces to a billiard in a polyhedral angle in $\R^n$. Global finiteness of the number of reflections in such billiards was proved by Ya.~G.~Sinai in the late 1970s; he also obtained an estimate depending only on the geometry of the angle; see \cite{Sinai78}.

Minkowski billiards generalize ``classical'' billiards. The difference is that the billiard table lies not in $\R^n$ but in an arbitrary finite-dimensional normed space with a smooth strictly convex, not necessarily symmetric, norm. In a normed space, the law ``the angle of incidence equals the angle of reflection'' makes no sense. It is replaced by a variational formulation: a billiard trajectory realizes an local extremum of the arc-length functional, with free points on the boundary of the billiard domain. Formal definitions are given in~\S\ref{sec:prelim}; see also~\cite{GT02}.

We do not rule out trajectories that ``hit a corner'' but regard them defined only up to the moment of such a hit, see definitions in \S\ref{sec:prelim}.

S.~V.~Ivanov and R.~V.~Barinov \cite{BI26}, \cite{IB26} previously proved local finiteness of reflections in semi-dispersing Minkowski billiards, i.e., billiards in the complement of a collection of convex sets. The results of \cite{BI26} used in this paper are stated in~\S\ref{sec:leg}.

The goal of this paper is to generalize the above mentioned results   on the number of reflections in polyhedral cones to Minkowski billiards. The main result is the following (Theorem~(\ref{t:finiteness})).

\begin{theorem}
Every trajectory of a Minkowski billiard in a polyhedral cone has only finitely many reflection points.
\end{theorem}

The proof uses a dimension reduction. We factor out the common line of the facets of the cone.  We then introduce the momentum of a trajectory, which is a natural integral of motion. For each momentum level, a new norm on the quotient space is constructed so that the projection of the original trajectory is a billiard trajectory in it; see~\S\ref{sec:pulse}.

We also show that global finiteness cannot be improved: we give an example of a Minkowski billiard with no uniform bound on the number of reflections; see~\S\ref{sec:ex}.

\textbf{Acknowledgments.} The author is sincerely grateful to his advisor, Professor S.~V.~Ivanov, for valuable ideas which greatly shaped this work.

\section{Preliminaries}\label{sec:prelim}

In this section we recall the basic definitions and state the results needed below; for details and proofs, see \cite{BI26}.

\subsection{Minkowski billiards}
\begin{definition}\label{d:minkowki norm}
A \textit{Minkowski space} is a finite-dimensional real vector space $X$ endowed with a \textit{Minkowski norm}, that is, a continuous function $\Norm\co X\to\R$ satisfying the following conditions:
\begin{enumerate}
\item
Positive definiteness: $\|x\|>0$ for all $x\in X\setminus\{0\}$.
\item
Positive homogeneity: $\|tx\|=t\|x\|$ for all $x\in X$ and $t\ge 0$.
\item
Subadditivity: $\|x+y\|\le \|x\|+\|y\|$ for all $x,y\in X$.
\item
Strict convexity, in the sense that the set $\{x\in X: \|x\|=1\}$ (called the \textit{indicatrix} or the \textit{unit sphere} of the norm) contains no line segments.
\item
$C^1$-smoothness on $X\setminus\{0\}$.
\end{enumerate}
\end{definition}

\begin{remark}
The following holds.
\begin{enumerate}
\item For positively homogeneous functions, the subadditivity is equivalent to convexity.
\item Strict convexity is equivalent to the subadditivity condition being strict for linearly independent vectors $x,y$.
\end{enumerate}
\end{remark}

In what follows we refer to Minkowski norms simply as norms; they are always assumed to be $C^1$-smooth. 
In the literature, Minkowski norms are usually assumed
to be $C^\infty$-smooth and quadratically strictly convex,
but for the purposes of this paper the above conditions are sufficient. Everywhere in the sequel the word ``smooth'' means $C^1$.

A \textit{Minkowski billiard} is determined by a Minkowski space $X=(X^n,\Norm)$ and a domain $\Omega\subset X$, called the \textit{billiard table}. We denote by $\ov\Omega$ the closure of $\Omega$, and by $\pd\Omega$ its boundary. A point $p\in\pd\Omega$ is called \textit{smooth} if, in a neighborhood of $p$, the boundary $\pd\Omega$ is a smooth hypersurface. For a smooth point $p\in\pd\Omega$, we denote by $T_p\pd\Omega$ the hyperplane in $X$ containing $0$ and parallel to the tangent hyperplane to $\pd\Omega$ at $p$.

A billiard trajectory can be informally described as a path of a particle that moves uniformly inside $\Omega$ and when it hits a smooth point $p\in\pd\Omega$, it reflects from the boundary according to the following variational principle: if $q_1,q_2$ are sufficiently close points of the trajectory before and after the reflection, then $p$ is a critical point of the function
\be\label{e:local length}
  x \mapsto \|x-q_1\| + \|q_2-x\|, \qquad x\in \pd\Omega.
\ee
We assume that trajectories are parametrized with unit speed. The formal details are given in the following definition.

\begin{definition}\label{d:trajectory}
The \textit{trajectory} of the Minkowski billiard $(X,\Norm,\Omega)$ is a continuous curve $\ga\co I\to\ov\Omega$, parametrized by an interval $I\subset\R$, such that for any interior point $t\in I$ there exists $\delta>0$ such that, for all $t_1\in(t-\delta,t)$ and $t_2\in(t,t+\delta)$ at least one of the following conditions holds:
\begin{enumerate}
\item[(1)]
The restriction $\ga|_{[t_1,t_2]}$ parametrizes a segment of an affine line in $X$ with unit speed (i.e.\ $\|\dot\ga\|=1$).
\item[(2)]
The point $p=\ga(t)$ is a smooth boundary point of $\pd\Omega$, the restrictions $\ga|_{[t_1,t]}$ and $\ga|_{[t,t_2]}$ are unit-speed parametrizations of the segments $[q_1,p]$ and $[p,q_2]$ where $q_1=\ga(t_1)$ and $q_2=\ga(t_2)$, and $p$ is a critical point of \eqref{e:local length} on $\pd\Omega$.
\end{enumerate}
The values of $t$ for which only the second condition holds are called \textit{reflection moments}, and the corresponding points $\ga(t)$ are called \textit{reflection points}. The \textit{length} of a trajectory $\ga$ is the length of the interval~$I$.

\end{definition}

We emphasize that the conditions in this definition apply only the interior points of the interval $I$; we do not impose any conditions on the behavior of the trajectory near its endpoints. One may assume that the trajectories are parametrized by open intervals and the endpoint values are defined by continuity if necessary.

\begin{definition}
A Minkowski billiard is called \textit{semi-dispersing} if the billiard table has the form
\be\label{e:walls}
 \Omega = X \setminus \bigcup_{i=1}^m W_i,
\ee
where $W_1,\dots,W_m\subset X$ are closed convex sets with nonempty interiors and smooth boundaries. The sets $W_i$ are referred to as the \textit{walls}.
\end{definition}

\subsection{The Legendre transform and the reflection law}\label{sec:leg}
\begin{definition}
The \textit{Legendre transform} of the norm $\Norm$ is the positively 1-homogeneous map $\LL=\LL_{\Norm}\co X\to X^*$ defined by
\be\label{e:legendre}
\LL(x) = d_x\big(\tfrac12\Norm^2\big) = \|x\|\, d_x(\Norm), \qquad x\in X.
\ee
The \textit{dual Legendre transform} of $\Norm$ is  also the positively 1-homogeneous map $\LL^*\co X^*\to X$ defined by the analogous formula, with respect to the dual norm $\Norm^*$:
\be\label{e:colegendre}
\LL^*(p) = d_p\big(\tfrac12(\Norm^*)^2\big) = \|p\|^*\, d_p(\Norm^*), \qquad p\in X^*.
\ee
Strictly speaking, $\LL(p)$ is an element of the double dual space $X^{**}$, but this space is canonically identified with $X$.
\end{definition}
If $\Norm$ is a Minkowski norm, one easily checks that the dual norm $\Norm^*$ is also a Minkowski norm. If $\|v\| = 1$, then the kernel of the functional $\LL(v)$ is the hyperplane through the origin parallel to the tangent to the unit sphere at $v$. Similarly, if $\|p\|^* = 1$, then $\ker\LL^*(p)$ is the analogous hyperplane with respect to $p$.

The following classical theorem is an important property of the Legendre transform.
\begin{theorem}\label{t:legendre inv}
Let $(X,\Norm)$ be a Minkowski space. Then
\begin{enumerate}
\item $\LL$ is a homeomorphism between the unit spheres of $\Norm$ and $\Norm^*$.
\item $\LL$ is an involution in the sense that $\LL^*=\LL^{-1}$.
\end{enumerate}
\end{theorem}

The following lemma reformulates the reflection law (\ref{e:local length}) of Minkowski billiards. In the sequel this reformulation will be used instead of the definition.

\begin{lemma}\label{l:reflection law}
Let $(X,\Norm,\Omega)$ be a Minkowski billiard, $p\in\ov\Omega$ a smooth boundary point, and $\ga$ a two-segment polygonal chain in $\ov\Omega$ with a break at~$p$, parametrized with unit speed. Let $v_-$ and $v_+$ be the velocity vectors of $\ga$ before and after the break.

Then $\ga$ is a billiard trajectory if and only if
\be\label{e:reflection}
 (\LL(v_+)-\LL(v_-))|_{T_p\pd\Omega} = 0.
\ee
\end{lemma}

The next statement asserts the uniqueness of the direction of a trajectory upon reflection off a wall.

\begin{lemma}\label{l:reflection exists}
Let $X=(X,\Norm)$ be a Minkowski space, let $H\subset X$ be a hyperplane containing~0, and let $v_-\in X\setminus H$ with $\|v_-\|=1$. Then there exists a unique vector $v_+\in X$ such that $v_+\ne v_-$, $\|v_+\|=1$, and
\be\label{e:reflection in H}
(\LL(v_+)-\LL(v_-))|_H = 0,
\ee
where $\LL$ is the Legendre transform of $\Norm$. Moreover, $v_+$ and $v_-$ lie in different open half-spaces with respect to~$H$.
\end{lemma}

By Lemmas \eqref{l:reflection law} and \eqref{l:reflection exists}, a billiard trajectory extends until it reaches a nonsmooth boundary point of the billiard table. More precisely,
\begin{cor}\label{cor:extension}
In a semi-dispersing Minkowski billiard, any maximal forward trajectory having a starting point either has infinite length or is a finite polygonal chain ending at a non-smooth boundary point.
\end{cor}

The main result of the paper \cite{BI26} is the following theorem, which will be used below.

\begin{theorem}\label{t:time}
In a semi-dispersing Minkowski billiard, the set of reflection moments along any billiard trajectory of finite length is finite.
\end{theorem}

\section{Billiards in polyhedral cones}

\subsection{Polyhedral cones}\label{sec:conus}

\begin{definition}\label{d:conus}
Let $\{H_i\}_{i=1}^N$ be a finite collection of hyperplanes through the origin of~$X$. For each $H_i$, fix a \textit{unit conormal} $\eta_i\in X^*$, that is a covector satisfying
\be\label{e:conormal}
\|\eta_i\|^* = 1,\quad \eta_i(h) = 0,\; \forall h \in H_i, \quad \eta_i(x)\neq 0,\; \forall x \in X\setminus H_i.
\ee
Define a \textit{positive (open) half-space} ${H_i^+:=\{x\co\, \eta_i(x) > 0\}}$; the hyperplane $H_i$ is its boundary.

A \textit{polyhedral cone} $\Omega$ with facets $\{H_i\}_{i=1}^N$ and conormals $\{\eta_i\}_{i=1}^N$ is the intersection of the corresponding closed positive half-spaces:
$$\Omega = \bigcap\limits^N_{i=1} \overline{H_i^+} = \{x\in X\co\,\eta_i(x)\ge 0,\,\forall i = 1,\dots,N\}.$$
A polyhedral cone $\Omega$ is called \textit{nondegenerate} if $\bigcap H_i = \{0\}$.

The collection $\{H_i,\eta_i\}_{i=1}^N$ is called \textit{minimal} if no hyperplane can be dropped without changing the cone:

$$
\nexists \,j\co\;\Omega = \bigcap\limits^N_{i=1} H_i^+ = \bigcap\limits^N_{i\neq j} H_i^+.
$$

Hereafter all collections defining polyhedral cones are assumed minimal.

It is easily to see that the term ``cone'' is justified, since the intersection of finitely many half-spaces is a cone with vertex at 0. Note also that every such cone is convex.
\end{definition}

\begin{example} The following examples illustrate the definition \ref{d:conus}.
\begin{enumerate}
\item $X=\R^2$. Let $H_1,\, H_2,\, H_3$ be three distinct lines through the origin. Then, regardless of the choice of positive half-planes, the collection $\{H_1,\, H_2,\, H_3\}$ is not minimal.
\item $X=\R^3$. Let $H_1,\, H_2$ be two distinct planes containing the line $z=0$. Then the cone with facets $H_1$ and $H_2$ is degenerate, regardless of the choice of conormals.
\item $X=\R^3$. Any cone whose facets form the lateral faces of a polyhedral pyramid is nondegenerate.
\end{enumerate}
\end{example}

\begin{remark}\label{r:conormal}
A collection of hyperplanes and conormals $\{H_i, \eta_i\}$ defines a nondegenerate cone if and only if the linear span of $\{\eta_i\}$ is the whole dual space~$X^*$.
\end{remark}
\begin{proof}
All the conormals $\eta_i$ lie in a proper subspace $N\subset X^*$ if and only if there exists a nonzero vector $v^{**}\in X^{**}$ such that $v^{**}(\eta_i)=0$ for all $i$ ($v^{**}$ is one of the conormals to $N$). Identifying $v^{**}\in X^{**}$ with a vector $v\in X$, we get $\eta_i(v)=0$; hence $v$ lies in all the hyperplanes $H_i$.
\end{proof}

Since half-spaces are convex and their boundaries are smooth hyperplanes, a Minkowski billiard whose table is a polyhedral cone is semi-dispersing; hence Theorem~\ref{t:time} applies to it.

\subsection{Reduction of the cone}\label{sec:reduction}
We fix some notation and conventions used in this and the following sections.

The vectors $v_-$ and $v_+$ conventionally denote the velocities before and after a given reflection. For definiteness, we take the velocity of a trajectory at a reflection time is $v_-$. The open unit ball and the unit sphere of $\Norm$ are denoted by $B$ and $S$ respectively, and those of the dual norm $\Norm^*$ by $B^*$ and $S^*$.

Let $\Gamma\subset\R$ be the maximal domain of definition of a trajectory $\ga$. In what follows we assume $\Gamma = (0,+\infty)$: if a trajectory does not extend to $+\infty$, then, by Theorem~\ref{t:time}, it has finitely many reflection points, and this case is of no interest. If a trajectory extends to infinity in both directions, it splits into two trajectories one of which corresponds to the reversed time and the opposite norm, defined by $\|x\|^-=\|{-x}\|$. Hence it suffices to consider trajectories with an initial point.

\begin{theorem}\label{t:number}
Let $\Omega\subset X$ be a nondegenerate polyhedral cone with facets $\{H_i\}_{i=1}^N$ and conormals $\{\eta_i\}_{i=1}^N$. Then any trajectory $\ga(t)$ of the Minkowski billiard $(X^n,\Norm,\Omega)$ eventually reflects off at most $N-1$ facets of $\Omega$:

$$
\exists\, T\in\R\co\;\{\text{reflection points of }\ga|_{(T,+\infty)}\}\subset \bigcup_{i\in I} H_i,\quad |I|\le N-1.
$$

\end{theorem}
\begin{proof}
Suppose the contrary; then there exists a trajectory $\ga$ such that for every $T\in\R$ there is a sequence of $k\ge N$ reflection times $t_1,\dots,t_k$, all greater than $T$, whose reflection points visit all $N$ facets, i.e.,
\begin{equation}\label{e:protivnoe}
\{\ga(t_1),\,\ga(t_2),\,\dots,\,\ga(t_k)\}\cap H_i \neq\emptyset,\quad \forall i \in \{1,\,2,\dots,N\}.
\end{equation}

Let $v_-$ and $v_+$ be the velocities of the trajectory before and after a reflection off the facet $H_j$. Then, by the reflection law \eqref{e:reflection}, for some nonzero $\alpha\in\R$ we have
\be\label{e:alpha}
\LL(v_+) = \LL(v_-) + \alpha\eta_j,
\ee
where $\alpha$ is chosen so that the covector $\LL(v_-)+\alpha\eta_j$ lies in $S^*$. By convexity of the dual sphere, the choice of $\alpha$ is unique.

Note that $\alpha>0$: indeed,

$$
1 = \LL(v_+)(v_+) = \LL(v_-)(v_+) + \alpha\eta_j(v_+),
$$

where the first summand on the right-hand side is strictly less than 1; hence the second one is positive. Moreover, $\eta_j(v_+)>0$, since $v_+$ points into the cone, hence lies in $H_j^+$.

Let $v_0$ be the velocity of the trajectory before the first reflection, and let $v_1,\,v_2,\,\dots,\,v_j,\,\dots$ be the velocities after the 1st, 2nd, \dots, $j$-th reflections respectively; let $H_{r_j}$ be the facet at which the $j$-th reflection occurs. From equality \eqref{e:alpha},

$$
\LL(v_{i}) = \LL(v_0)+\sum_{j=1}^{i-1}\alpha_j\eta_{r_j},
$$

where $\eta_{r_j}$ is the conormal corresponding to the facet $H_{r_j}$, and $\alpha_j$ is the corresponding constant defined by \eqref{e:alpha}.

Choose a unit vector $e$ lying in the interior of the cone $\Omega$; it exists by nondegeneracy. Then, from the preceding equality, the fact that the velocity vectors are unit, the positivity of $\alpha$, and the choice of $e$, we get

$$
2\ge \LL(v_{i+1})(e) - \LL(v_0)(e) = \sum^i_{j=1}\alpha_j\eta_{r_j}(e)\ge 0.
$$

Hence, for some $m\le N$ with $\eta_m(e)$ minimal among the conormals corresponding to the facets,

$$
2\ge \sum^i_{j=1}\alpha_j\eta_{r_j}(e) \ge \eta_{m}(e)\sum^i_{j=1}\alpha_j\ge 0.
$$

Since the partial sums are bounded and $\alpha_j>0$, the series $\sum^\infty_{j=1}\alpha_j$ converges. From the equality
$$
\|\LL(v_{i+1}) - \LL(v_i)\|^* = \alpha_i\|\eta_{r_j}\|^* = \alpha_i,
$$
it follows that the series $\sum^\infty_{i=1} \|\LL(v_{i+1}) - \LL(v_i)\|^*$ converges; hence its tails are small:
\begin{equation}\label{e:ryad}
\forall\,\varepsilon>0\; \exists\, m\co\;\sum^\infty_{i=m} \|\LL(v_{i+1}) - \LL(v_i)\|^*\,<\,\varepsilon.
\end{equation}

Given $\delta>0$, define the \textit{$\delta$-almost equator in the direction $p\in X^*$} as the set

$$
\delta E(p)\coloneqq\{q\in S^*\co (q+\varepsilon p)\in S^* \text{ for some } 0<\varepsilon\le\delta\};
$$

informally, these are the points $q\in S^*$ from which the ray in the direction $p$ meets $S^*$ again within distance $\delta$.

Define also $E_0(p)$, the \textit{equator in the direction $p\in X^*$}: the set of points $q\in S^*$ such that the direction $p$ is tangent to $S^*$ at $q$. Smoothness, compactness and convexity of $S^*$ imply
\be\label{e:equator}
E_0(p) = \bigcap_{\delta>0} \overline{\delta E(p)}.
\ee

Let $v_-$ and $v_+$ be the velocities of the trajectory before and after a reflection off the facet $H_k$, such that $\|\LL(v_+)-\LL(v_-)\|^* \le \delta$. Then equality \eqref{e:alpha} implies that $\LL(v_-)$ lies in $\delta E(\eta_k)$.

Define the function $l(\delta)$ as the minimal length of a curve visiting all the closed almost equators $\overline{\delta E(\eta_i)}$; the minimum exists by compactness of the sphere. By definition, $l(\delta)$ does not decrease as $\delta$ decreases.

Note that there exists $\delta_0$ such that $l(\delta_0)>0$;  otherwise, for every $\delta>0$, the sets $\ov{\delta E(\eta_i)}$ share a point $p\in S^*$; by \eqref{e:equator}, so do all the equators $E_0(\eta_i)$. Then all the conormals $\eta_i$ are tangent to $S^*$ at $p$, hence lie in one hyperplane of $X^*$, contradicting Remark~\ref{r:conormal}.

\begin{figure}[h!]
   \centering
   \begin{tikzpicture}[
    tangent/.style={
        decoration={
            markings,
            mark=
                at position #1
                with
                {
                    \coordinate (tangent point-\pgfkeysvalueof{/pgf/decoration/mark info/sequence number}) at (0pt,0pt);
                    \coordinate (tangent unit vector-\pgfkeysvalueof{/pgf/decoration/mark info/sequence number}) at (1,0pt);
                    \coordinate (tangent orthogonal unit vector-\pgfkeysvalueof{/pgf/decoration/mark info/sequence number}) at (0pt,1);
                }
        },
        postaction=decorate
    },
    use tangent/.style={
        shift=(tangent point-#1),
        x=(tangent unit vector-#1),
        y=(tangent orthogonal unit vector-#1)
    },
    use tangent/.default=1
]

    \draw[tangent=0, tangent=0.1, tangent=0.5, tangent=0.6] (0,-4) to[out = 180, in = -90] (-3,0) to [out=90,in=180] (0,2) to [out=0,in=90] (3,0) to [out=-90,in=0] (0,-4);
    \node at (2.,-3.5){$S^*$};


    \draw[name path = E2f] (tangent point-1) to [out=180,in=180] (tangent point-3);
    \node at (-1.,-2.25){$E_0(\eta_2)$};
    \draw[loosely dashed,name path = E2b] (tangent point-1) to [out=0,in=-90] (2.5,0.2) to [out = 90, in = 0] (tangent point-3);

    \draw[name path = E3f] (tangent point-2) to [out=-45,in=-35] (tangent point-4);
    \node at (1.5,-2.75){$E_0(\eta_3)$};
    \draw[loosely dashed,name path = E3b] (tangent point-2) to [out=135,in=-120] (-1.75,-0.1) to [out = 60, in = 155] (tangent point-4);

    \draw[name path = E1f] (-3,0) arc [start angle=-180, end angle=0, x radius=3, y radius=1] node[pos=0.55,above]{$E_0(\eta_1)$};
    
    \draw[name path = E1b, loosely dashed] (-3,0) arc [start angle=180, end angle=0, x radius=3, y radius=1];
    
    \path[name intersections={of=E1f and E2f, by=I}];
    \fill (I) circle (1pt);
    
    \path[name intersections={of=E2f and E3f, by=K}];
    \fill (K) circle (1pt);

    \path[name intersections={of=E1f and E3f, by=L}];
    \fill (L) circle (1pt);

    \draw[use tangent = 3, ->, >=latex] (0,0) -- ({-1.5*cos(-10)},{1.5*sin(10)}) node[ above right] {$\eta_2$};    
    \draw[use tangent = 1, shift = (I), ->, >=latex] (0,0) -- ({1.5*cos(-10)},{1.5*sin(-10)}) node[pos=0.8, below left] {$\eta_2$};
    \draw[use tangent = 1, shift = (K), ->, >=latex] (0,0) -- ({1.5*cos(-10)},{1.5*sin(-10)}) node[pos=0.8, below ] {$\eta_2$};

    \draw[use tangent = 2, ->, >=latex, shift = (K)] (0,0)--(-1.5,0) node[left, pos = 0.8]{$\eta_3$};
    \draw[use tangent = 2, ->, >=latex, shift = (L)] (0,0)--(-1.5,0) node[left]{$\eta_3$};    
    \draw[use tangent = 2, ->, >=latex, shift = (tangent point-4)] (0,0)--(-1.5,0) node[right]{$\eta_3$};

    \draw[->, >=latex] (-3,0)--(-3,1.5) node[right]{$\eta_1$};
    \draw[->, >=latex] (L)--($(L)+(0,1.5)$) node[below left]{$\eta_1$};
    \draw[->, >=latex] (I)--($(I)+(0,1.5)$) node[below left]{$\eta_1$};
    
\end{tikzpicture}
     \caption{Arrangement of the equators $E_0(\eta_i)$ on the dual sphere}
 \end{figure}
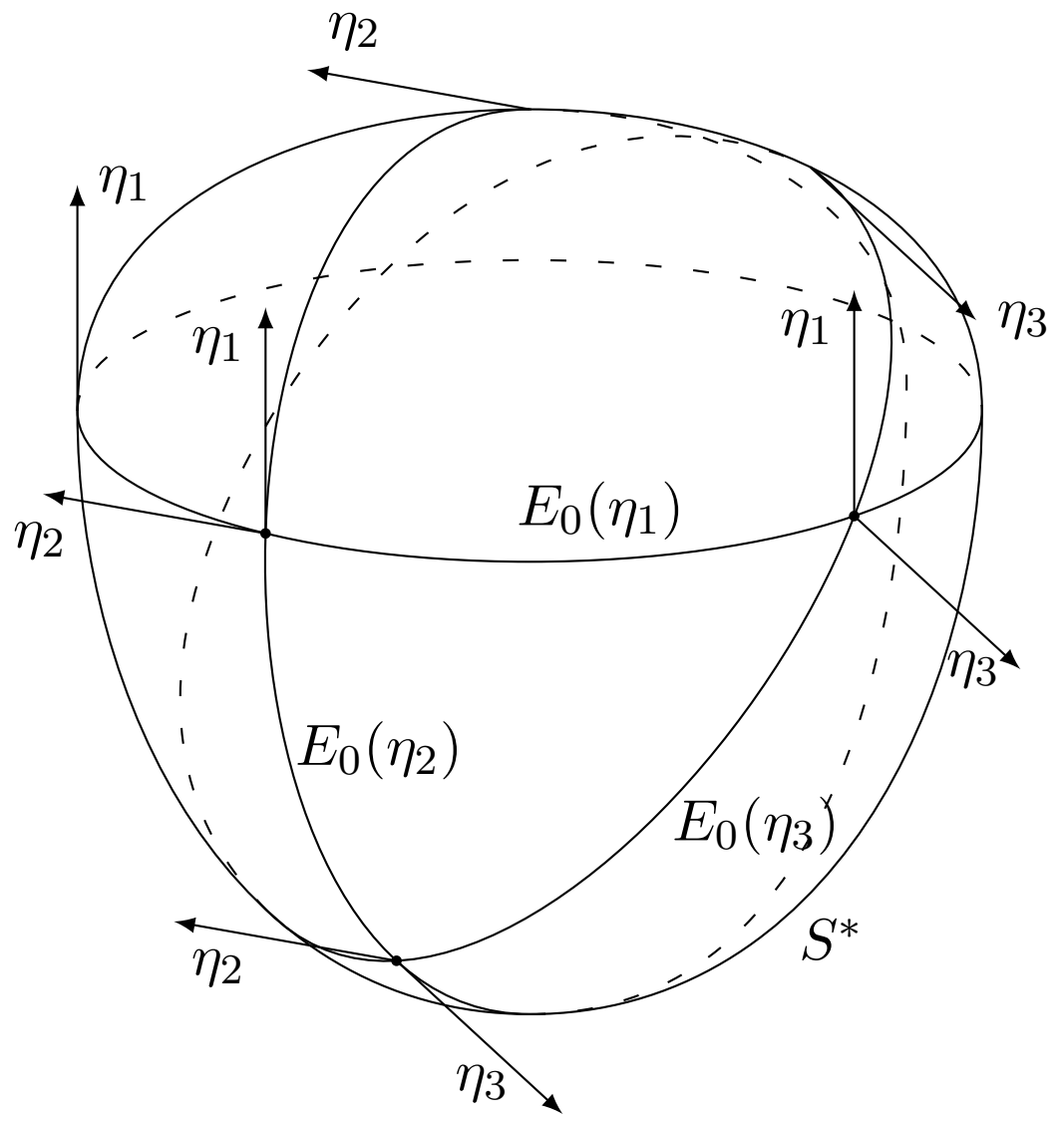

Since the series in \eqref{e:ryad} converges, there exists $M$ such that
\be\label{e:dlina lomannoy}
\sum^\infty_{i=M} \|\LL(v_{i+1}) - \LL(v_i)\|^*\,<\min(\delta_0, l(\delta_0)).
\ee

By assumption \eqref{e:protivnoe}, a sufficiently long sequence of reflections of the trajectory, starting from the $M$-th one, visits all $N$ facets of the cone. Then the sequence $\{\LL(v_i),\; i\ge M\}$ visits all $N$ almost equators $\delta_0E(\eta_i)$; hence the left-hand side of \eqref{e:dlina lomannoy} is at least the length of the curve connecting the vectors of this sequence as points on the sphere, and therefore at least $l(\delta_0)$ — a contradiction:

$$
l(\delta_0)\le\sum^\infty_{i=M} \|\LL(v_{i+1}) - \LL(v_i)\|^*\,<\min(\delta_0, l(\delta_0)).
$$

\end{proof}

\begin{remark}
The above proof generalizes the argument from the paper of Ya.~G.~Sinai \cite{Sinai78}.
\end{remark}

\subsection{Momentum of a trajectory and dimension reduction}\label{sec:pulse}

\qquad Let $E\subset X$ be a subspace of $X$. There exists a canonical projection $\pi\co X\to X/E$ such that $\pi(x)=[x]$, where $[x]$ is the coset of $x$. The operator $\pi^*\co (X/E)^*\to X^*$, dual to the projection, is given by

$$
\pi^*\varphi(x) = \varphi([x]),\quad \forall\,\varphi\in (X/E)^*;
$$

that is, $\pi^*\varphi$ is a functional on $X$ that lifts a functional from the quotient space by applying it to cosets. It is well known that $\pi^*$ is a canonical isomorphism between $(X/E)^*$ and the annihilator $E^\perp$.

For a set of vectors $M\subset X$, denote its image under the projection by $[M]=\pi(M)$.

\begin{theorem}\label{t:finiteness}
    Let $\Omega\subset X^n$ be a polyhedral cone with facets $\{H_i\}_{i=1}^{N}$ and conormals $\{\eta_i\}_{i=1}^{N}$. Then every trajectory $\ga(t)$ of the Minkowski billiard $(X^n,\Norm,\Omega)$ has finitely many reflections:
    
    $$
    \exists\,T\co\;\ga(T,+\infty)\cap H_i = \emptyset,\quad \forall\,i = 1,\dots,N.
    $$

\end{theorem}

By Theorem~\ref{t:number}, it suffices to prove Theorem~\ref{t:finiteness} for degenerate cones only. Indeed, after a finite time the trajectory will reflect only from facets sharing a common line, and by Theorem~\ref{t:time} only finitely many reflections can occur in finite time.

The proof proceeds by induction on dimension, which is lowered by factoring out the line $E$ common to all facets. We first introduce an invariant of the trajectory and use it to construct a new Minkowski norm on the quotient for which the induction hypothesis holds.

\begin{definition}
    The \textit{momentum of a trajectory} $\ga(t)$ \textit{with respect to a vector} $e$ is the function $P_{\ga,e}(t)$ given by
    \begin{equation*}
        P_{\ga,e}(t)\coloneqq \LL(\ga'(t))(e).
    \end{equation*}
\end{definition}

\begin{lemma}\label{l:pulse}
    Let a billiard trajectory $\ga(t)$ in a polyhedral cone reflect off facets sharing a common line $E$, and let $e$ be a unit vector generating $E$. Then the momentum of $\ga$ with respect to $e$ is constant:
    \begin{equation*}
        P_{\ga,e}(t) = \LL(\ga'(t))(e) = C(\ga,e),
    \end{equation*}
    where $C(\ga,e)$ is a constant depending only on the trajectory and the sign of vector $e$.
\end{lemma}
\begin{proof}
    By the reflection law \eqref{e:reflection}, for any billiard trajectory $\ga(t)$ and any $t_1,\,t_2\in\Gamma$ we have
    
    $$
    \LL(\ga'(t))(e)\big|_{t=t_1} = \LL(\ga'(t))(e)\big|_{t=t_2},
    $$

    because before and after a reflection the covectors $\LL(v_-)$ and $\LL(v_+)$ agree on vectors of the reflecting facet, and $e$ lies in every facet.
\end{proof}
\begin{remark}
     The momentum $P_{\ga,e}$ is an invariant of the trajectory: by Noether's theorem, it corresponds to the symmetry of translations along $e$.
\end{remark}

From now on, $E$ denotes the line common to all facets, and $e$ a unit vector generating it.

\paragraph{Construction of the norm on $\bm{X/E}$.}\label{hyperplanes}
The family of affine hyperplanes $\{L_\alpha\}_{\alpha\in[-1,1]}\subset X^*$ is defined by

$$
 L_\alpha\coloneqq \{p\in X^*\co\, p(e) = \alpha\}. 
$$
All the $L_\alpha$ are parallel to one another and, in particular, to $L_0$, which passes through the origin. By definition \eqref{e:conormal}, $L_0$ is the $(n-1)$-dimensional subspace containing the linear span of the conormals $\{\eta_1,\dots,\eta_N\}$. Note that $L_0$ is the annihilator of $E$, so $L_0 = E^\perp \cong (X/E)^*$, where the second isomorphism is the canonical identification.

\begin{definition}\label{d:Ca}
Fix an arbitrary covector $a\in (B^*\cap L_C)$, that is a covector such that $\|a\|^*<1$ and $a(e) = C$. Define the \textit{dual norm} $\Norm^*_{C,a}$ on $(X/E)^*$ by
\begin{equation}\label{e:Ca}
    \|p\|^*_{C,a} = \|\pi\8 p+a\|^*,\quad p\in(X/E)^*.
\end{equation}
Geometrically, the unit sphere of $\Norm^*_{C,a}$ is obtained by intersecting $S^*$ with $L_C$, shifting the result by $a$, and pulling it back under $\pi\8$ (see Fig.~\ref{pic:Norm_Ca}).

The \textit{norm} $\Norm_{C,a}$ on $X/E$ is the norm corresponding to the dual norm $\Norm^*_{C,a}$. The Legendre transforms corresponding to the new norms are denoted $\LL_{C,a}$ and $\LL_{C,a}^*$.
\end{definition}
The norms constructed are indeed Minkowski norms, since so is the original norm $\Norm$.

\begin{lemma}\label{l:equivalent}
    Let $\LL(v)(e) = C$, and let $[\ov v]$ be the unit vector of $\Norm_{C,a}$ in the direction of $[v]$. Then
    
    $$
    \LL(v)-a = \pi\8 \LL_{C,a}([\ov v]),\;[\ov v] = \frac{[v]}{\|[v]\|_{C,a}}.
    $$

\end{lemma}
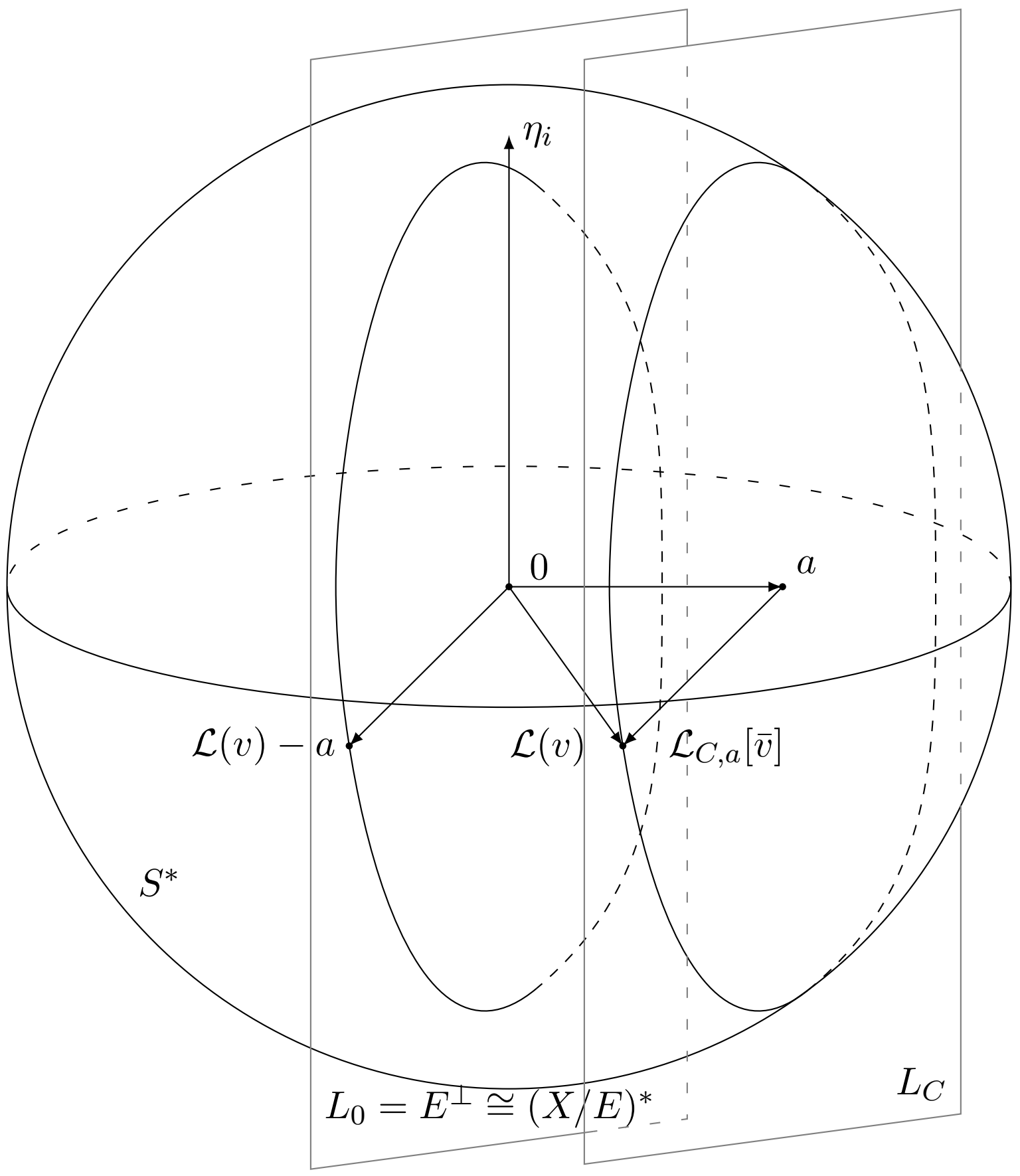
\begin{figure}[h]
  \centering

  \begin{tikzpicture}[
    tangent/.style={
        decoration={
            markings,
            mark=
                at position #1
                with
                {
                    \coordinate (tangent point-\pgfkeysvalueof{/pgf/decoration/mark info/sequence number}) at (0pt,0pt);
                    \coordinate (tangent unit vector-\pgfkeysvalueof{/pgf/decoration/mark info/sequence number}) at (1,0pt);
                    \coordinate (tangent orthogonal unit vector-\pgfkeysvalueof{/pgf/decoration/mark info/sequence number}) at (0pt,1);
                }
        },
        postaction=decorate
    },
    use tangent/.style={
        shift=(tangent point-#1),
        x=(tangent unit vector-#1),
        y=(tangent orthogonal unit vector-#1)
    },
    use tangent/.default=1
]

    \draw (0,0) circle (5);
    \node[above right] at ({5*cos(-140)},{5*sin(-140)}){$S^*$};
    \coordinate (a) at (2.725,0);
    \coordinate (O) at (0,0);
    \begin{scope}
    \clip (-5, -2) rectangle (5, 0);    
    \draw (0,0) ellipse (5 and 1.2);
    \end{scope}
        \begin{scope}
    \clip (-5, 2) rectangle (5, 0);    
    \draw[loosely dashed] (0,0) ellipse (5 and 1.2);
    \end{scope}

    \draw[tangent=0.65] (3,4) to [out=140,in=90] (1,0) to [out = 270, in = -140] (3,-4);
    \draw[dashed] (3,4) to [out = -40, in = 90] (4.25,0) to [out = -90, in = 40] (3,-4);
    \fill[] (tangent point-1) circle (1pt);
    \draw[->, >=latex] (a) -- (tangent point-1) node[pos=0.8,below right] {\(\LL_{C,a}[\bar{v}]\)};
    \draw[->,>=latex] (O) -- (tangent point-1) node[pos = 0.8, below left] {$\LL(v)$};    
    
    \draw[tangent=0.65] (0.275,4) to [out=140,in=90] (-1.725,0) to [out = 270, in = -140] (0.275,-4);
    \draw[dashed] (0.275,4) to [out = -40, in = 90] (1.525,0) to [out = -90, in = 40] (0.275,-4);
    \fill[] (tangent point-1) circle (1pt);
    \draw[->, >=latex] (0,0) -- (tangent point-1) node[left] {\(\LL(v)-a\)};
    
    \fill (0,0) circle (1pt);
    \node at (0.3, 0.2) {\(0\)};

    \fill (a) circle (1pt) node[above right]{$a$};
    
    \draw[->, >=latex] (0,0) -- (0, 4.5) node[right] {\(\eta_i\)};
    \draw[->,>=latex] (O) -- (a);

    \draw[gray] (1.775, 5.5) -- (1.775, 5.75) -- (-1.975, 5.25) -- (-1.975, -5.8)  node[black, above right, pos = 0.975]{\(L_0 = E^\perp \cong (X/E)^*\)} -- (0.8375, -5.425);
    \draw[loosely dashed,gray, on background layer] (1.775, 5.5) -- (1.775, -5.3) -- (0.8375, -5.425);
        
    \draw[gray] (4.5, 2.25) -- (4.5, 5.75) -- (0.75, 5.25) -- (0.75, -5.75) -- (4.5, -5.25) -- (4.5, -2.25) node[black, above left, pos = 0]{$L_C$};
    \draw[loosely dashed,gray, on background layer] (4.5,2.25) -- (4.5, -2.25); 
\end{tikzpicture}
\caption{Relation between the norms $\Norm\8$ and $\Norm\8_{C,a}$}
\label{pic:Norm_Ca}
\end{figure}
\begin{proof}
    Consider two elements of the bidual $(X/E)^{**}$:
    \begin{enumerate}
        \item $[\ov v]^{**} = \LL\8_{C,a}(\LL_{C,a}[\ov v])$, the normalized projection of the vector $v$ transformed twice by the Legendre transform with respect to the new norm;
        \item $[v^{**}] = [\LL\8(\LL( v))]$, the projection onto $(X/E)^{**}$ of the original vector $v$ transformed twice by the Legendre transform with respect to the original norm.
    \end{enumerate}

    By the involution property of $\LL$ and the properties of the projection,
    \begin{equation}\label{be}
    \begin{aligned}
    \pi\8\ker ([\ov v]^{**}) &=\pi\8\ker ([v]^{**}) = \{p\in L_0\co\;p(v)=0 \},\\
    \pi\8\ker[v^{**}] &= \ker v^{**}\big|_{L_0}= \{p\in L_0\co\;p(v)=0 \},
    \end{aligned}
    \end{equation}
    i.e., the images under $\pi\8$ of the kernels of the covectors $[\ov v]^{**}$ and $[v^{**}]$ coincide. Note that $\pi\8\ker ([\ov v]^{**})$ is a linear $2$-codimensional subspace, parallel to the tangent hyperplane to $\pi\8S\8_{C,a}$ at the point $\pi\8\LL_{C,a}[\ov v]$, where $S\8_{C,a}$ is the dual sphere of the norm $\Norm\8_{C,a}$; and $\pi\8\ker[v^{**}]$ is the intersection of $L_0$ with the hyperplane parallel to the tangent hyperplane to $S\8$ at $\LL(v)$.

    The smooth strictly convex sets $(L_C\cap S\8) -a$ and $\pi\8S\8_{C,a}$ coincide (by construction); by the above, $\LL(v)-a$ and $\pi\8\LL_{C,a}[\ov v]$ are points of this set corresponding to the same tangent direction. Similarly to the equality of kernels \eqref{be}, one obtains the equality of the positive half-spaces; hence $\LL(v)-a$ and $\pi\8\LL_{C,a}[\ov v]$ coincide.
\end{proof}

\begin{cor}\label{cor:bil}
    Let $[\ov\ga(t)]$ be the natural reparametrization of $[\ga(t)]$, where $\ga(t)$ is a billiard trajectory, and let $P_{\ga,e} = C$. Then the trajectory $[\ov\ga(t)]$ is a billiard trajectory in the space $(X/E,\, \Omega/E,\,\Norm_{C,a})$.
\end{cor}
\begin{proof}
    Let $[\ov v_-]$ and $[\ov v_+]$ be the velocities of $[\ov\ga(t)]$ before and after a reflection off $H/E$; by Lemma~\ref{l:equivalent},
    
    $$
    (\pi\8)^{-1}\LL(v_-)-a = \LL_{C,a}([\ov v_-]),\qquad (\pi\8)^{-1}\LL(v_+)-a = \LL_{C,a}([\ov v_+]).
    $$

    By the reflection law \eqref{e:reflection}, the trajectory $[\ov\ga(t)]$ is a billiard trajectory, since $\ga$ is billiard and
    
    $$
    \big(\LL_{C,a}([\ov v_+])-\LL_{C,a}([\ov v_-])\big)\big|_{H/E} = \big((\pi\8)^{-1}\LL(v_+) - (\pi\8)^{-1}\LL(v_-)\big)\big|_{H/E} = 0.
    $$
\end{proof}

\paragraph{Proof of Theorem \ref{t:finiteness}.}

\textit{Base case.} For $n=1$, it is obvious.


Suppose the theorem holds for all $i\le n-1$; we show that it holds for $n$.

\textit{Induction step.} If $\Omega$ is nondegenerate, discard facets one by one, by Theorem~\ref{t:number}, until the cone becomes degenerate. Let $E$ be the common line of the remaining facets. Pass to the quotient by $E$. Consider image of $\ga(t)$ in $X/E$, reparametrized naturally in the norm $\Norm_{C,a}$ with $C = P_{\ga,e}$. It has the same number of reflections as $\ga$, and by Corollary~\ref{cor:bil} it is a billiard trajectory in $X/E$; hence, by the induction hypothesis, it has finitely many reflections.
\qed

\begin{remark}
    The construction of the new norm was the key step of the proof: the standard ``projection'' of the original norm onto the quotient not only fails to preserve billiard trajectories, but may fail to be a convex norm at all.
\end{remark}

\section{No uniform bound on the number of reflections}\label{sec:ex}

This section presents an example of a Minkowski billiard with no uniform bound on the number of reflections of its trajectories: for every $n\in\mathbb{N}$ there exists a trajectory with more than $n$ reflections. This is an essential difference from the Euclidean case, since for any semi-dispersing Euclidean billiard there is a bound on the number of reflections of every trajectory, depending only on the metric and the geometry of the billiard domain (see \cite{BFK98}).

Consider the space $\R^3$ with the standard Euclidean metric; for convenience we denote its dual by $\R^{3*}$ and identify it with $\R^3$ in the standard way. The billiard domain is the wedge (dihedral angle) 
$$
K\coloneqq \{(x,y,z)\in \R^3\co\; z>|y|\}
$$
 \label{wedge} with facets $H_1\coloneqq\{y=z,\,z\ge0\}$ and $H_2\coloneqq\{y=-z,\,z\ge0\}$, and edge $E\coloneqq\{y=z=0\}$.

The norm is defined by its dual sphere $S^*$, whose construction can be informally described as follows: the standard sphere is truncated, in a neighborhood of its intersection with the $x$-axis, by an ellipsoid that is flat in the direction of the $x$-axis and stretched along the $z$-axis; the resulting truncation is then truncated again by an even flatter ellipsoid stretched even more along the $z$-axis, and so on; passing to the limit and smoothing the resulting surface, we obtain $S^*$. A more formal description of $S^*$ is given below.

\subsection{Construction of the norm}

\begin{definition}
    Consider the sequence of ellipsoids $\{E_i\subset\R^{3*},\, i\in\mathbb{N}_0\}$ given by the equations
\be
\frac{x^2}{a_i}+\frac{y^2}{b_i}+\frac{z^2}{c_i} = 1,
\ee
where the coefficients satisfy the following properties:

\begin{enumerate}\label{d:ellipsoids}
    \item The ellipsoid $E_0$ is the standard unit sphere, i.e., $a_0=b_0=c_0=1$.
    \item The ellipsoids are flattened along the $x$-axis, i.e., there is a constant $A\in (0,1)$ such that $a_i\searrow A$.
    \item The ellipsoids are stretched along the $y$- and $z$-axes, i.e., ${c_i}\nearrow \infty$ and ${b_i}\nearrow \infty$ as $i\to\infty$.
    \item The ellipsoids are stretched along the $z$-axis more strongly than along the $y$-axis, i.e., $\frac{c_i}{b_i}\nearrow \infty$ as $i\to\infty$.
    \item In the half-space $\{x\ge0\}$, the intersections $E_{i-1}\cap E_i$ and $E_i\cap E_{i+1}$ are separated by some plane perpendicular to the $x$-axis, i.e.,
\be\label{5cond}\max\limits_{E_{i-1}\cap E_i} x < \min\limits_{E_{i}\cap E_{i+1}} |x|.\ee
\end{enumerate}
\end{definition}
\begin{figure}[h]
    \centering
    \begin{tikzpicture}[scale = 0.95]

\newcommand{\ellArc}[4]{%
  \path (#1); \pgfgetlastxy{\ax}{\ay}%
  \path (#2); \pgfgetlastxy{\bx}{\by}%
  \pgfmathsetmacro{\axcm}{\ax/28.45274}%
  \pgfmathsetmacro{\aycm}{\ay/28.45274}%
  \pgfmathsetmacro{\bxcm}{\bx/28.45274}%
  \pgfmathsetmacro{\bycm}{\by/28.45274}%
  \pgfmathsetmacro{\angA}{atan2(\aycm/#4, \axcm/#3)}%
  \pgfmathsetmacro{\angB}{atan2(\bycm/#4, \bxcm/#3)}%
  \draw[thick] (#1) arc[start angle=\angA, end angle=\angB,
                        x radius=#3, y radius=#4]
}

\newcommand{\ellArcRev}[4]{%
  \path (#1); \pgfgetlastxy{\ax}{\ay}%
  \path (#2); \pgfgetlastxy{\bx}{\by}%
  \pgfmathsetmacro{\axcm}{\ax/28.45274}%
  \pgfmathsetmacro{\aycm}{\ay/28.45274}%
  \pgfmathsetmacro{\bxcm}{\bx/28.45274}%
  \pgfmathsetmacro{\bycm}{\by/28.45274}%
  \pgfmathsetmacro{\angA}{atan2(\aycm/#4, \axcm/#3)}%
  \pgfmathsetmacro{\angB}{atan2(\bycm/#4, \bxcm/#3)}%
  \pgfmathsetmacro{\angB}{\angB < \angA ? \angB + 360 : \angB}%
  \pgfmathsetmacro{\deltaB}{\angB - \angA}%
  \pgfmathsetmacro{\angA}{\deltaB < 180 ? \angA - 360 : \angA}%
  \draw[ thick] (#1) arc[start angle=\angA, end angle=\angB,
                        x radius=#3, y radius=#4]
}

\clip (-6.1,-6.1) rectangle (6.1,6.1);

\draw[dashed, gray, name path=E0](0,0) circle (5);
\draw[dashed, gray, name path=E1](0,0) ellipse (4 and 6);
\draw[dashed, gray, name path=E2](0,0) ellipse (3.5 and 15);

\path[name intersections={of=E0 and E1, name=i, total=\t}];
\path[name intersections={of=E1 and E2, name=j, total=\u}];

\coordinate (P1) at (i-1);
\coordinate (P4) at (i-4);
\coordinate (P2) at (j-1);
\coordinate (P3) at (j-4);
\coordinate (P5) at (-5,0);
\coordinate (P6) at (0,-5);

\ellArc{P1}{P2}{4}{6};
\ellArc{P3}{P4}{4}{6};

\ellArc{P2}{P3}{3.5}{15};
\ellArcRev{P1}{P4}{5}{5};
\foreach \p in {P1,P2,P3,P4}{\fill (\p) circle (1.2pt);}
\node[below right] at (0.0,-5){$E_0$};
\node[above left] at (3.2,-3.6){$E_1$};
\node[below left] at (3.5,-0.5){$E_2$};
\node[above right] at (0,0){$0$};
\fill (0,0) circle (1.2pt);
\draw[->,>=latex] (-0.5,0)--(5.5,0) node[below]{$x$};
\draw[->,>=latex] (0,-0.5)--(0,5.5) node[right]{$z$};
\end{tikzpicture}
\caption{The section of the ellipsoids $E_0,\,E_1$, and $E_2$ by the plane $xOz$}
    \label{pic:ell}
\end{figure}
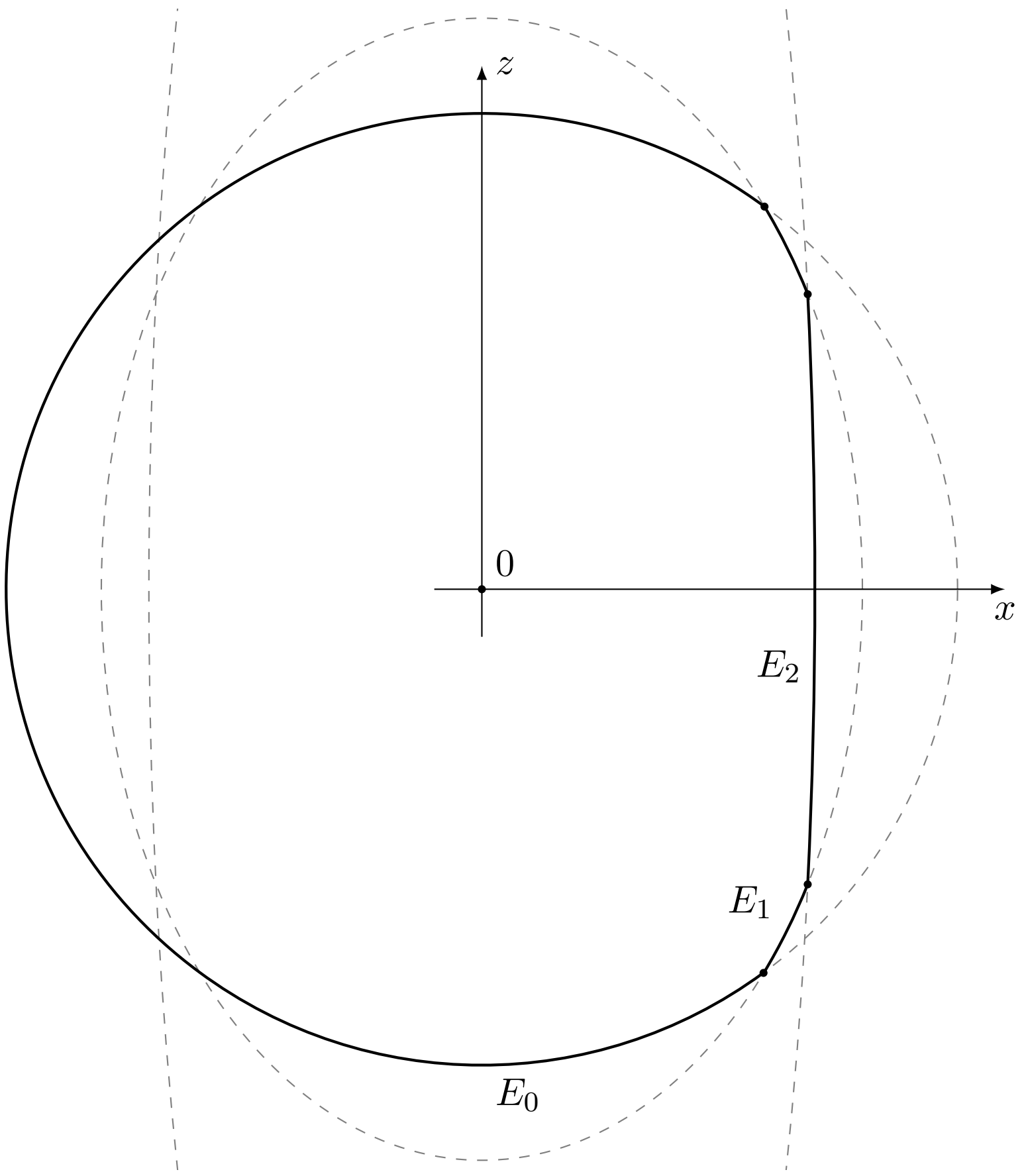
\begin{remark}
    Geometrically, the ellipsoids are chosen as follows. The standard sphere $E_0$, in a neighborhood of the point $(1,0,0)$, is cut by the ellipsoid $E_1$, which is flat in the direction of the $x$-axis and stretched along the other two. Next, for an ellipsoid $E_2$, stretched along the $y$- and $z$-axes more strongly than $E_1$, the coefficient $a_2$ is chosen so that $E_2$ intersects $E_1$ in a neighborhood of the point $(a_1,0,0)$ so small that property~5 of Definition~\ref{d:ellipsoids} holds. Such an $a_2$ exists: for $a_2=a_1$ the ellipsoids $E_1$ and $E_2$ are tangent, so a slight decrease of $a_2$ makes inequality \eqref{5cond} hold. The subsequent ellipsoids are chosen in the same way. (See Fig.~\ref{pic:ell}.)
\end{remark}

Such a sequence $\{E_i\}$ can be given explicitly.

\begin{example}
The sequence of ellipsoids $\{E_i,\, i\in\mathbb{N}_0\}$ given by the equations
\be
\frac{x^2}{a_i}+\frac{y^2}{b_i}+\frac{z^2}{c_i} = 1,\quad \mbox{where}\quad a_i = 1-\left(\sum\limits_{k=1}^{2i} \frac{1}{4^k}\right)=\frac{2}{3}+\frac{1}{3\cdot16^i},\;b_i=\left(\frac{9}{4}\right)^i,\;c_i=4^i
\ee
satisfies all the conditions of Definition~\ref{d:ellipsoids}.
\end{example}
The verification of this claim is purely arithmetic; being tedious, it is omitted.

\begin{definition}
Let $m_i$ and $M_i$ denote the minimum and maximum of $|x|$ on $E_{i-1}\cap E_i$, respectively, and let $\Delta_i = M_i - m_i$. By inequality \eqref{5cond}, $m_i<M_i<m_{i+1}$ and $\Delta_i>0$ for all $i$. The \textit{dual sphere} $S^*$ is defined by the following rules:
\begin{enumerate}
    \item For $x<m_1-\sfrac{\Delta_1}{10}$, the sphere $S^*$ coincides with $E_0$, the standard unit sphere.
    \item For $M_i+\sfrac{\min(\Delta_i, \Delta_{i+1})}{10}<x<m_{i+1}-\sfrac{\min(\Delta_i, \Delta_{i+1})}{10}$, the sphere $S^*$ coincides with the ellipsoid $E_i$.
    \item For $m_i-\sfrac{\min(\Delta_i, \Delta_{i-1})}{10}\le x\le M_i+\sfrac{\min(\Delta_i, \Delta_{i+1})}{10}$, the sphere $S^*$ is the convex smoothing joining $E_{i-1}$ and $E_i$ while respecting the previous two rules.
    \item For $x=A$, the sphere $S^*$ consists of the single point $(A,0,0)$.
\end{enumerate}
Note that the resulting sphere is smooth and convex. \textit{The norm $\Norm_S$ on $\R^3$} is defined so that $S^*$ is its dual sphere.

In Fig.~\ref{pic:ell}, the black line shows the section of the unsmoothed sphere $S^*$ for the ellipsoids $\{E_0,\,E_1,\,E_2\}$.
\end{definition}

\subsection{Proof of the absence of a uniform bound}

\begin{theorem}
    In the Minkowski billiard $(\R^3,\, K,\,\Norm_S)$, for every natural number $n$ there exists a trajectory with at least $n$ reflections.
\end{theorem}
\begin{proof}
    Let $e$ be the unit vector (in the norm $\Norm_S$) generating the edge of the wedge $K$, i.e., the $x$-axis. Consider the planes $\{L_\alpha\subset\R^{3*}\}$ associated with the norm $\Norm_S$ as in \ref{hyperplanes}; note that they are perpendicular to the $x$-axis (in the standard metric on $\R^{3*}$), since they are parallel to $L_1$, which is tangent to $S^*$ at the point $(A,\,0,\,0)$ by the construction of $\Norm_S$. Introduce coordinates $(y,\,z)$ on each $L_\alpha$ as the natural projection of the standard coordinates on $\R^{3*}$ along the $x$-axis.

    Choose $\varepsilon_i\in(0,\,1)$ so that the section of the sphere $S^*$ by the plane $L_{\varepsilon_i}$ is the ellipse~$\mathcal{E}_i^*$ given by
    
    $$
    \frac{y^2}{b_i}+\frac{z^2}{c_i} = \operatorname{const},
    $$

    where $b_i$ and $c_i$ are the coefficients from Definition~\ref{d:ellipsoids}. Define the norm $\Norm^*_{\varepsilon_i,O_i}$ on $(\R^3/E)^*$ as in Definition~\ref{d:Ca}, where $O_i$ is the intersection point of $L_{\varepsilon_i}$ with the $x$-axis. Note that this is a Euclidean norm, since its unit sphere is the ellipse $\mathcal{E}_i^*$ centered at the origin. Then the unit sphere of the corresponding norm $\Norm_{\varepsilon_i,O_i}$ is the ellipse $\mathcal{E}_i$ given by
    
    $$
    \frac{y^2}{c_i}+\frac{z^2}{b_i} = \operatorname{const}.
    $$

    Consider the billiard $(\R^3/E,\,K/E,\, \Norm_{\varepsilon_i,O_i})$: this is a Euclidean billiard in the right angle~$\{z\ge|y|\}$. Apply a linear transformation of the plane taking $\mathcal{E}_i$ to the standard circle; the billiard trajectory (in the norm $\Norm_{\varepsilon_i,O_i}$) becomes a standard one, that is one with equal angles of incidence and reflection. The right angle $\{z\ge|y|\}$ becomes an acute angle of magnitude $\alpha_i$, where $\alpha_i\to 0$ as $i\to\infty$, since $\frac{c_i}{b_i}\to\infty$. The following elementary lemma gives the maximal number of reflections of a trajectory in such billiards; its proof is sketched below.

    \begin{lemma}\label{flat}
Any trajectory of a standard two-dimensional billiard in an angle of size $\alpha$ is reflected at most $N(\alpha)$ times, where
\be\label{Nalpha}
N(\alpha)=
\begin{cases}
\left[\dfrac{\pi}{\alpha}\right]+1, & \text{if } \dfrac{\pi}{\alpha}\notin\mathbb{Z} ,\\[6pt]
\left[\dfrac{\pi}{\alpha}\right], & \text{if } \dfrac{\pi}{\alpha}\in\mathbb{Z},
\end{cases}
\ee
and the bound is sharp.
\end{lemma}

For any natural $i$, consider a trajectory $\ga_i$ with $P_{\ga_i,e} = \varepsilon_i$ in the billiard $(\R^3,\, K,\,\Norm_S)$ such that its projection in the billiard $(\R^3/E,\,K/E,\, \Norm_{\varepsilon_i,O_i})$ has exactly $N(\alpha_i)$ reflections; note that the trajectory $\ga_i$ has the same number of them. Since $\alpha_i\underset{i \to \infty}{\to} 0$, we have $N(\alpha_i)\underset{i \to \infty}{\to} \infty$, and therefore

$$
\text{the number of reflections}(\ga_i)\underset{i \to \infty}{\to} \infty.
$$

\end{proof}

\textit{Proof of Lemma \ref{flat}.} A detailed proof is given in [\cite{GaZ}, \S15, A.]; here we sketch the idea.

Denote the angle by $Q_0$ and its sides by $a_0$ and $b_0$. Fix a trajectory $\ga_0$; without loss of generality, its first reflection occurs off the side $a_0$. Reflect the angle $Q_0$ and the tail of $\ga_0$ after this reflection across the line $a_0$. The result is a new angle $Q_1$ sharing the side $a_0$ with $Q_0$; its other side $b_1$ is the mirror image of $b_0$. The pre-reflection part of $\ga_0$ together with the segment of $\ga_1$ from $a_0$ to $b_1$ forms one straight line. Repeat this procedure: reflect across $b_1$, then across $a_2$, and so on, each time across the side hit by the trajectory. Stop when the sum $Q_1+Q_2+\dots$ (not counting $Q_0$) exceeds $\pi$; this happens after exactly $N(\alpha)$ reflections. (If the reflections stop earlier, we stop too.)

Note that the trajectory ``straightens out'' step by step, so the $(N(\alpha)+1)$-th reflection does not occur. Moreover, along any line intersecting the angle $Q_0$, a billiard trajectory is uniquely recovered by the inverse sequence of reflections. Then choose a line intersecting all the $Q_i$; it corresponds to a trajectory with exactly $N(\alpha)$ reflections.
\qed

\end{document}